\documentclass{amsart}
\usepackage{mathrsfs,comment,setspace}

\newtheorem{theorem}{Theorem}[section]

\theoremstyle{definition}

\theoremstyle{remark}
\newtheorem{remark}[theorem]{Remark}

\numberwithin{equation}{section}

\begin{document}

\title[Student's $t$-distribution]{A simple proof of the characteristic function of Student's $t$-distribution}

\author{Robert E. Gaunt}
\address{School of Mathematics, The University of Manchester, Manchester M13 9PL, UK}
\thanks{The author is supported by a Dame Kathleen Ollerenshaw Research Fellowship.}





\begin{abstract}
This note presents a simple proof of the characteristic function of Student's $t$-distribution. The method of proof, which involves finding a differential equation satisfied by the characteristic function, is applicable to many other distributions.
\end{abstract}

\maketitle

\vspace{-10mm}

\section{Introduction}Let $X\sim t_{\nu}$ follow Student's $t$-distribution with $\nu>0$ degrees of freedom and density function
\begin{equation} \label{fancy} f_X(x)=\frac{\Gamma(\frac{\nu+1}{2})}{\sqrt{\pi \nu}\Gamma(\frac{\nu}{2})}\bigg(1+\frac{x^2}{\nu}\bigg)^{-\frac{1}{2}(\nu+1)}, \quad x\in\mathbb{R}.
\end{equation}
Over the years, the characteristic function of Student's $t$-distribution has received much interest in the statistics literature.  Starting in 1956 from the work of Fisher and Healy \cite{fh56}, complicated series form solutions that are only valid for odd, even or fractional degrees of freedom, resulting from complicated and sometimes convoluted proofs, were given by Ifram \cite{i70}, Pestana \cite{p77}, Mitra \cite{m78} and Sutradha \cite{s86}; see also Johnson, Kotz and Balakrishnan \cite{jkb95} for further details regarding the expressions of Mitra.  Finally, independently in 1995, Hurst \cite{h95} and Joarder \cite{j95} gave an elegant closed-form formula that is valid for all $\nu>0$: 
\begin{equation}\label{result}\mathbb{E}[\mathrm{e}^{\mathrm{i}tX}]=\frac{K_{\nu/2}(\sqrt{\nu} |t|)(\sqrt{\nu} |t|)^{\nu/2}}{\Gamma(\frac{\nu}{2})2^{\nu/2-1}}, \quad t\in\mathbb{R},
\end{equation}
where $K_\nu(x)$ is a modified Bessel function of the second kind (see Appendix \ref{appa} for a definition and basic properties that are needed in this note).  Hurst \cite{h95} arrived at this formula as a limiting case of the symmetric generalized hyperbolic distribution characteristic function, whilst Joarder \cite{j95} used a standard integral representation of $K_\nu(x)$ to obtain a very efficient proof.  Since then, Joarder and Alam \cite{j95b} obtained the characteristic function of the elliptic $t$-distribution, Joarder and Ali \cite{j96} derived the characteristic function of the multivariate $t$-distribution, and Dreier and Kotz \cite{dk02} used the theory of positive definite densities to obtain a new integral representation of the classic Student's $t$-distribution characteristic function.

Motivated by the historical interest in this problem, the purpose of this short note is to complement the existing literature with a new simple and direct derivation of the characteristic function of Student's $t$-distribution.  As part of the proof we find an ordinary differential equation (ODE) that the characteristic function must satisfy, which is similar to the classical modified Bessel function differential equation.  This gives a transparent explanation as to why the characteristic function is given in terms of the modified Bessel function $K_\nu(x)$.  Also, as will be elaborated on in Remark \ref{rem4}, the proof serves as a useful exposition of a technique that could be used to provide simple derivations of characteristic functions of many other distributions.







\section{Proof}



Let $X\sim t_{\nu}$ and denote $\phi_X(t)=\mathbb{E}[\mathrm{e}^{\mathrm{i}tX}]$.  As the distribution of $X$ is symmetric about $0$, $\phi_X(-t)=\phi_X(t)$ for all $t\in\mathbb{R}$.  It therefore suffices to consider the case $t>0$.  The rest of the proof consists of two parts.  We first prove (\ref{result}) for $\nu>2$ and then deduce that the formula must also be valid for all $\nu>0$.  Suppose $\nu>2$.  As $(\nu+x^2)f_X'(x)+(\nu+1)xf_X(x)=0$, an integration by parts gives that
\begin{align*}0&=\int_{-\infty}^\infty \big[(\nu+x^2)f_X'(x)+(\nu+1)xf_X(x)\big]\mathrm{e}^{\mathrm{i}tx}\,\mathrm{d}x\\
&=\int_{-\infty}^\infty \big[-\mathrm{i}t(\nu+x^2)+(\nu-1)x\big]\mathrm{e}^{\mathrm{i}tx}f_X(x)\,\mathrm{d}x \\
&=-\mathrm{i}t\mathbb{E}[X^2\mathrm{e}^{\mathrm{i}tX}]+(\nu-1)\mathbb{E}[X\mathrm{e}^{\mathrm{i}tX}]-\mathrm{i}\nu t\mathbb{E}[\mathrm{e}^{\mathrm{i}tX}],
\end{align*}
where all integrals exist because $\nu>2$.  Since $\phi_X'(t)=\mathrm{i}\mathbb{E}[X\mathrm{e}^{\mathrm{i}tX}]$ and $\phi_X''(t)=-\mathbb{E}[X^2\mathrm{e}^{\mathrm{i}tX}]$, it follows that $\phi_X(t)$ satisfies the ODE
\begin{equation}\label{plo}t\phi_X''(t)-(\nu-1)\phi_X'(t)-\nu t\phi_X(t)=0.
\end{equation}
It follows from (\ref{realfeel}) that the general solution to (\ref{plo}) is  $\phi_X(t)=C_1 t^{\nu/2}I_{\nu/2}(\sqrt{\nu} t)+C_2 t^{\nu/2}K_{\nu/2}(\sqrt{\nu} t)$, where $C_1$ and $C_2$ are arbitrary constants.  This is because direct differentiation shows that if $f(x)$ satisfies $x^2f''(x)+xf'(x)-(x^2+(\nu/2)^2)f(x)=0$ then $g(x)=x^{\nu/2}f(\sqrt{\nu}x)$ satisfies $xg''(x)-(\nu-1)g'(x)-\nu xg(x)=0$. Being a characteristic function, $\phi_X(t)$ must satisfy the conditions $|\phi_X(t)|<\infty$ for all $t\in\mathbb{R}$ and $\phi_X(0)=1$, and we can use the limiting forms (\ref{roots}) and (\ref{Ktend0}) to determine $C_1$ and $C_2$ accordingly.  This yields formula (\ref{result}) for the case $t>0$, $\nu>2$.  


Now suppose $\nu>0$.  Let $Y\sim t_{\nu+2}$, whose characteristic function $\phi_Y(t)=\mathbb{E}[\mathrm{e}^{\mathrm{i}tY}]$ we have already found in the first part of the proof.  Then  
\begin{align*}\phi_X(t)=\mathbb{E}[\mathrm{e}^{\mathrm{i}tX}]=-\frac{\mathrm{i}}{t}\sqrt{\frac{\nu}{\nu+2}}\mathbb{E}[Y\mathrm{e}^{\mathrm{i}t\sqrt{\frac{\nu}{\nu+2}}Y}]=-\frac{1}{t}\sqrt{\frac{\nu}{\nu+2}}\phi_{Y}'\Big(\sqrt{\frac{\nu}{\nu+2}}t\Big),
\end{align*}
where the second equality can be obtained from an integration by parts followed by a rescaling of the integration variable and a simplification of a constant using the standard formula $\Gamma(u+1)=u\Gamma(u)$.
On using (\ref{ddbk}) to compute $\phi_Y'(\sqrt{\frac{\nu}{\nu+2}}t)$ we obtain formula (\ref{result})  (after again using that $\Gamma(u+1)=u\Gamma(u)$) for the case $t>0$, $\nu>0$, which completes the proof. \hfill $\Box$

\begin{remark}\label{rem4}It is evident that the proof technique can be used to obtain characteristic functions for many other continuous random variables. One simply needs to find an ODE with polynomial coefficients satisfied by the probability density function, and then use integration by parts to find an ODE satisfied by the characteristic function.  (This can be done provided the distribution satisfies certain integrability conditions (for example, we needed to take $\nu>2$ in this step of our proof).)  If the ODE is tractable then it possible to deduce a formula for the characteristic function.  In some cases (as in the Student's $t$ case)  we may first need to prove the result for a restricted range of parameters before extending to the full range of validity.  An example in which this can be done using an argument similar to the one we used for Student's $t$-distribution is the $F$-distribution with $d_1$ and $d_2$ degrees of freedom with density proportional to $x^{d_1/2-1}(1+d_1x/d_2)^{-(d_1+d_2)/2}$, $x>0$ (for the characteristic function see Phillips \cite{p82}).  
\end{remark}

\appendix

\section{Modified Bessel functions}\label{appa}

The following basic properties of modified Bessel functions can be found in Olver et al$.$ \cite{olver}, Chapter 10.  The general solution of the modified Bessel differential equation
\begin{equation} \label{realfeel} x^2 f''(x) + xf'(x) - (x^2 +\nu^2)f(x) =0\end{equation}
is given by $f(x)=C_1I_{\nu} (x) +C_2K_{\nu} (x)$, where the modified Bessel functions of the first kind $I_\nu(x)$ and second kind $K_\nu(x)$ are defined, for $\nu\in\mathbb{R}$ and $x>0$, by
\[I_{\nu} (x) =  \sum_{k=0}^{\infty} \frac{(\frac{1}{2}x)^{\nu+2k}}{\Gamma(\nu +k+1) k!} \quad \text{and} \quad K_\nu(x)=\frac{1}{2}\int_0^\infty u^{\nu-1}\mathrm{e}^{-x(u+u^{-1})/2}\,\mathrm{d}u.
\]
For $\nu>0$, the modified Bessel functions have the following properties:
\begin{eqnarray} \label{roots} I_{\nu} (x) &\sim& \frac{\mathrm{e}^x}{\sqrt{2\pi x}}, \quad x \rightarrow \infty, \\
\label{Ktend0}K_{\nu} (x) &\sim&  2^{\nu -1} \Gamma (\nu) x^{-\nu},  \quad x \downarrow 0,  \\
 \label{ddbk}\frac{\mathrm{d}}{\mathrm{d}x}\big(x^\nu K_\nu(x)\big)&=&-x^{\nu} K_{\nu-1}(x).
\end{eqnarray}

\bibliographystyle{amsplain}

\end{document}